\documentclass[12pt]{amsart}
\usepackage{amsmath,amsfonts,amscd,amssymb}
\newtheorem{theorem}{Theorem}

\newtheorem{corollary}{Corollary}
\theoremstyle{remark}\newtheorem{remark}{Remark}
\theoremstyle{remark}

\begin{document}
\title[Edge colorings and Hamiltonian cycles]
{The enumeration of edge colorings and Hamiltonian cycles
by means of symmetric tensors}

\author{Peter Zograf}
\address{Steklov Mathematical Institute, St.Petersburg, 191011 Russia}
\email{zograf@pdmi.ras.ru}

\thanks{Partially supported by the ARC grant DP0343028.} 
\subjclass[2000]{Primary 05C45; Secondary 05C15}
\keywords{Graph function, edge coloring, Hamiltonian cycle}

\begin{abstract}
Following Penrose, we introduce a family of graph functions defined 
in terms of contractions of certain products of symmetric tensors 
along the edges of a graph. Special cases of these functions 
enumerate edge colorings and cycles of
arbitrary length in graphs (in particular, Hamiltonian  cycles).
\end{abstract}

\maketitle

\section{Introduction}The number of Hamiltonian cycles of a graph is an
elusive invariant that is hard to compute.  One can find several formulas
in the literature that express the number of Hamiltonian cycles of a graph in
terms of its adjacency matrix. However, all these formulas are rather 
complicated (e.g., they involve either traces of powers \cite{L} or determinats 
and permanents \cite{G-J} of submatrices of the adjacency matrix) and
cannot be efficiently used for counting Hamiltonian cycles. The aim of this
note is to propose a different approach to the enumeration of 
Hamiltonian cycles in graphs based on tensor contractions. The adjacency 
matrix  enters into this construction somewhat implicitly by governing the 
the order in which the tensors are contracted. The same construction works
for the enumeration of edge colorings and cycles of arbitrary length.

The idea of constructing graph invariants by means of tensor contractions  
traces back to Penrose \cite{P}, who proposed a novel (though not quite
successful) approach to the 4-color problem.  Later this idea found a number 
of remarkable applications, e.g.  in the theory of Vassiliev knot invariants 
\cite{K, BN2} (some other applications are also mentioned in the
survey \cite{CDK}).
However, they all deal with 3-valent embedded graphs (like
3-valent planar maps in the 4-color problem or Feynman diagrams in the 
theory of Vassiliev knot invariants) and make use of antisymmetric tensors 
(the structure tensors of Lie algebras, to be precise). Here we adapt Penrose's
construction to arbitrary graphs and symmetric tensors and apply it to the
above mentioned enumeration problems in graphs.

{\bf Acknowledgements.} The author is grateful to the Institute of Mathematical 
Sciences at Stony Brook University for hospitality and support during 2002-03 
academic year. 

\section{Symmetric tensors and graph functions}\label{function}
Let $G$ be a finite graph (possibly with loops and multiple edges). 
The set of vertices 
of $G$ we denote by $V(G)=\{v_1, \dots ,v_n\}$,  where $n=|V(G)|$ 
is the total number of vertices, and
the set of edges of $G$ we denote by $E(G)$. For each vertex $v_i$
we denote by $d_i$ its degree (or valency), $i=1, \dots ,n$.
Then the number of edges of $G$ is given by 
$$|E(G)|=\frac{1}{2}\sum_{i=1}^n d_i.$$

Now let $\mathbb{F}$ be a field, and let $V\cong\mathbb{F}^r$ be a vector 
space of dimension $r$ over $\mathbb{F}$. Fix a symmetric bilinear form 
$B: V\otimes V\longrightarrow \mathbb{F}$. The graph $G$ 
together with the bilinear form $B$ define a multilinear form
\begin{equation}
B_G: V^{\otimes d_1}\otimes \dots \otimes V^{\otimes d_n}
\longrightarrow\mathbb{F},
\end{equation}
which is constructed as follows. At each vertex $v_i$ of $G$ we place  
$d_i$\nobreakdash-th tensor power $V^{\otimes d_i}$ of the vector space 
$V$, where the factors are labeled by the half-edges of $G$ incident to $v_i$. 
Each edge of $G$ defines a contraction of two copies of $V$ (corresponding 
to its two half-edges) by means of the bilinear
form $B$.  We obtain the multilinear form $B_G$ by performing 
such contractions over the set $E(G)$ of all edges of $G$.
Rigorously speaking, the multilinear form $B_G$ depends on the order of
half-edges at each vertex $v_i$, or, equvalently, on the order of factors in the
tensor power $V^{\otimes d_i}$. However, its restriction to $S^{d_1} V\otimes
\dots\otimes S^{d_n} V$, where $S^d V$ denotes the $d$\nobreakdash-th 
symmetric power of $V$, is defined uniquely. 

Now fix a sequence $\mathcal{A}=\{A_1,A_2, \dots \}$ of symmetric contravariant 
$d$\nobreakdash-valent tensors $A_d\in S^d V \subset V^{\otimes d}$. We treat  
the tensor product $A_{d_1}\otimes \dots \otimes A_{d_n}$ as an element of 
$V^{\otimes d_1}\otimes \dots \otimes V^{\otimes d_n}$ 
and consider the element
\begin{equation}
\mathcal{F}_{\mathcal{A},\, B}(G)
=B_G (A_{d_1}\otimes \dots \otimes A_{d_i}) \in \mathbb{F}.
\end{equation}
Roughly speaking, $\mathcal{F}_{\mathcal{A},\, B}(G)$ is obtained by 
placing
a copy of $A_d$ at each vertex of $G$ of degree $d$ and contracting 
$\otimes_{i=1}^n A_{d_i}$
using $B$ over $|E(G)|$ pairs of indices corresponding to the edges of $G$.
Thus, to each pair $\mathcal{A},\, B$, where $\mathcal{A}$ is a sequence of
symmetric $d$-tensors ($d=1,2, \dots $) and $B$ is a
symmetric bilinear form, we associate an $\mathbb{F}$-valued mapping 
$\mathcal{F}_{\mathcal{A},\, B}$ on the set of isomorphism 
classes of graphs, or an  
$\mathbb{F}$-valued {\em graph function} in the terminology of \cite{T}.

\section{Enumeration of edge colorings}

Given a graph $G$, an {\em r-edge coloring} of $G$ is a coloring of
edges of $G$ in $r$ colors such that at each vertex $v_i\in V(G)$ 
all $d_i$ edges
incident to it have different colors. Clearly, an $r$-edge coloring
exists only if $r\geq d_i$ for all $i=1,\dots, n=|V(G)|$ and 
only if $G$ contains no loops. 

Let us show that the number of $r$-edge colorings of $G$ can be
realized as a graph function $\mathcal{F}_{\mathcal{A},\, B}$
for some suitably chosen $\mathcal{A}$ and $B$.
Put $\mathbb{F}=\mathbb{R}$ 
and consider the usual Euclidean coordinates in $V=\mathbb{R}^r$. 
Take $B$ given in these coordinates by the identity 
$r\times r$-matrix $I_r$, and define the components of the tensors 
$A_d,\; d=1,2,\dots,$ by the formula
$$A_d^{i_1\dots i_d}=\begin{cases}1& \text{if $i_1,\dots, i_d\in\{1,\dots, r\}$ 
are pairwise distinct,}\\ 
0& \text{otherwise}. \end{cases}$$

\begin{theorem}
For $\mathcal{A}$ as above,
the value of the graph function $\mathcal{F}_{\mathcal{A},\, I_r}$
on any graph $G$ is equal to the number of $r$-edge colorings of $G$.
\label{color}\end{theorem}

\begin{proof}
In order to compute the value $\mathcal{F}_{\mathcal{A},\, I_r}(G)$
first we have to decide which products of components of $A_d$
contribute to it non-trivially. We interprete the indices $1,\dots ,r$ as
colors of the half-edges of $G$. A product of $n$ components  
$A_{d_1}^{i_1 \dots\, i_{d_1}} \dots A_{d_n}^{i_{m-d_n+1} \dots\, i_m}$ 
(where $n=|V(G)|$ is the number of vertices and 
$m=\sum_{j=1}^n d_j=2|E(G)|$ is twice the number of edges of $G$) 
makes a non-zero contribution to $\mathcal{F}_{\mathcal{A},\, I_r}(G)$ 
if and only if the colors agree on each edge of $G$ or, equivalently, 
if and only if for each edge the both indices that label two 
of its half-edges are the same (otherwise the bilinear form
$B=I_r$ vanishes). By definition, $A_{d}^{i_1 \dots\, i_d}\neq 0$
if and only if all the indices $i_1 \dots\, i_d\in\{1,\dots,r\}$ are distinct, 
so that the non-zero contributions are in one-to-one
correspondence with $r$-edge colorings of $G$. Since every such
contribution is equal to 1,  the number of  $r$-edge colorings of $G$
is exactly $\mathcal{F}_{\mathcal{A},\, I_3}(G)$.
\end{proof}

\begin{remark}The case of 3-edge colorings of 3-valent graphs (also
called {\em Tait colorings}), is of a special interest  because of its
relation to the 4-color problem. Namely, 
a planar 3-valent map is 4-colorable if and only if it admits a Tait 
coloring, cf. \cite{P, BN1, CDK}. Take $\mathcal{A} =\{A_3\}$
and $B=I_3$ in Theorem 1. Then for any 3-valent graph $G$ the value 
$\mathcal{F}_{A_3,\, I_3}(G)$ is equal to the number of 
Tait colorings of $G$. 
\end{remark}

\section{Enumeration of cycles}

Given a graph $G$, by a {\em (multi)cycle} we understand  a 2-valent
subgraph $C$ in $G$. Note that we do not require a cycle to be connected. 
We denote by $|C|$ the {\it length}
of the cycle $C$ (that is, the number of edges of $G$ that belong to $C$), 
and by $l(C)$ the number of connected componets of $C$. 
Clearly, the length of a cycle cannot exceed $n=|V(G)|$,
and cycles of length $n$ are spanning cycles in $G$.  
A {\em Hamiltonian cycle} is a {\em connected} cycle of length $n$.

The {\em type} of a cycle $C$ in $G$ is the partition 
$\lambda_C=[|C_1|,\dots,|C_l|]$
of the number $|C|$, where $C_1,\dots,C_l$ are the connected componets of
$C,\;l=l(C)$. The weight of partition $\lambda_C$ is $|\lambda_C|=|C|$, and
the length is $l(\lambda_C)=l(C)$. For each partition $\lambda$ we define
a graph function $N_\lambda$ by 
$$N_\lambda(G)=\#\{C\subset G|\lambda_C=\lambda\},$$
i.e., $N_\lambda(G)$ is the number of cycles
of type $\lambda$ in $G$. 
In particular, $N_{[n]}(G)$ is the number of Hamiltonian cycles in $G$.

For $k$ a positive integer, denote by $p_k(x_1,x_2,\dots)=x_1^k+x_2^k+\dots$ 
the $k$\nobreakdash-th power sum in variables $x_1,x_2,\dots$.
Given a partition $\lambda=[k_1,\dots,k_l]$, we define a homogeneous
symmetric function $p_\lambda$ of degree $|\lambda|=k_1+\dots+k_l$ 
by the formula
$$p_\lambda(x_1,x_2,\dots)=\prod_{i=1}^l\;p_{k_i}(x_1,x_2,\dots).$$

We want to show that under a special choice of $\mathcal{A}$ and $B$ 
the graph function $\mathcal{F}_{\mathcal{A},\, B}$ defined in Section 
\ref{function} counts the number of cycles of any given type in graphs. 
We take $\mathbb{F}=\mathbb{C}$ and consider the standard 
coordinates in $V=\mathbb{C}^r$. In these coordinates 
the bilinear form $B$
is given by the identity $r\times r$ matrix $I_r$.
We put 
$$A_1^i =\begin{cases} 0 & \text{if $i\neq r$},\\ 
t & \text{if $i=r$},\end{cases}$$
and define the tensors $A_d$ for $d\geq 2$ componentwise
by the formula
$$A_d^{i_1 \dots \,i_d}=\begin{cases} x_i & \text{if $(i_1 \dots \,i_d)$ is a 
permutation of $(i\,i\,r \dots r)$},\\ 
& \hspace{2in} i=1, \dots ,r-1,\\ 
t & \text{if $(i_1 \dots \,i_d)=(r\,\dots\, r)$},\\
0 & \text{otherwise},\end{cases}$$
where $x_1,\dots, x_{r-1}$ and $t$ are arbitrary complex numbers.
The main result of this section is the following

\begin{theorem}
For $\mathcal{A}=\{A_1,A_2,\dots\}$ as above, 
the value of the graph function
$\mathcal{F}_{\mathcal{A}, \, I_r}$ on any graph $G$ is given 
by the formula
\begin{equation}
\mathcal{F}_{\mathcal{A}, \, I_r}(G)=\sum_{|\lambda|\leq n}
t^{n-|\lambda|}p_\lambda(x_1,\dots,x_{r-1})\, N_\lambda(G),
\label{poly}\end{equation}
where the sum is taken over the set of all partitions $\lambda$
of weight $|\lambda|\leq n=|V(G)|$.
\label{main}\end{theorem}

\begin{proof} 
As in the proof of Teorem 1, we interprete the indices $1,\dots ,r$ as
colors of the half-edges of $G$. Similarly, a product of $n$ components  
$A_{d_1}^{i_1 \dots\, i_{d_1}} \dots A_{d_n}^{i_{m-d_n+1} \dots\, i_m}$ 
(where $n=|V(G)|$ and $m=\sum_{j=1}^n d_j=2|E(G)|$) 
contributes non-trivially to $\mathcal{F}_{\mathcal{A},\, I_r}(G)$ 
if and only if the colors agree on each edge of $G$ or, equivalently, 
if and only if for each edge the both indices that label two 
of its half-edges are the same. Thus, in this case 
the non-zero contributions are in one-to-one
correspondence with edge colorings of $G$ in $r$ colors 
with the following properties:

(i) an edge incident to a vertex of degree 1 has color $r$, and

(ii) at each vertex $v_j\in V(G)$ of degree $d_j\geq2$ 
two edges incident to it have  some  
color $i_j\in\{1, \dots ,r\}$, and the rest $d_i-2$ edges have color 
$r$ (if an edge makes a loop we count it twice).

The closure of the union of edges with colors $1, \dots, r-1$
is a cycle $C$ in $G$,  and every connected component $C_j$ of $C$
is colored in one of the colors $i_j\in\{1, \dots ,r-1\}$. 
The contribution to $\mathcal{F}_{\mathcal{A},\, I_r}(G)$ from this coloring
is $t^{n-|\lambda_C|}\prod_{j=1}^{l(C)} x_{i_j}^{|C_j|}$, where 
$\lambda_C=[|C_1|,\dots, |C_l|]$ is the partition associated with $C$
and $l=l(C)$ is the number of connected components of $C$. 
Therefore, the contribution from all possible colorings of the cycle 
$C$ is equal to 
$$t^{n-|\lambda_C|}\prod_{j=1}^{l(C)}\left(\sum_{i=1}^{r-1} x_i^{|C_j|}\right)
=t^{n-|\lambda_C|}p_{\lambda_C}(x_1,\dots, x_{r-1}),$$
and summig up the contributions from all cycles in $G$ we get the 
assertion of the theorem.
\end{proof}

\begin{corollary}
The graph function $\mathcal{F}_{\mathcal{A}, \, I_r}$,
depending on $x_1,\dots, x_{r-1}$ and $t$ as parameters,  determines
the numbers $N_\lambda(G)$ uniquely for any graph
$G$ with $n\leq r-1$ vertices.
\end{corollary}

\begin{proof}
By Theorem \ref{main},  the graph function 
$\mathcal{F}_{\mathcal{A}, \, I_r}$  with values in $\mathbb{C}$
factors through the ring
$\mathbb{C}[x_1,\dots, x_{r-1}]^{S_{r-1}}$ of 
symmetric polynomials  
in $r-1$ independent variables $x_1,\dots, x_{r-1}$.
It is well known that the polynomials $p_k(x_1,\dots, x_{r-1}),\;
k=1,\dots, n$, are algebraically independent in
$\mathbb{C}[x_1,\dots, x_{r-1}]^{S_{r-1}}$ 
provided $n\leq r-1$. Therefore, in this case the graph function
$\mathcal{F}_{\mathcal{A}, \, I_r}(G)$ determines the coefficients
$N_\lambda(G)$ in (\ref{poly}) uniquely.
\end{proof}

\begin{remark}
Since the coefficients $N_\lambda(G)$ in (\ref{poly}) are
{\em non-negative} integers, we can uniquely find them out when
$$r\geq 1+\underset{C\subset G}{\max}\; l(C),$$
where $l(C)$ is the number of  the connected componenets of $C$
and the maximum is taken over all cycles $C$ in $G$,
but we will not dwell on this here.
\label{remark}\end{remark}

Below are two special cases of Theorem \ref{main} of independent
interest.

\begin{corollary}
Put $r=2$, $x_1=1$ and $t=0$. Then for any graph $G$
the value $\mathcal{F}_{\mathcal{A}, \, I_2}(G)$ is
the number of spanning cycles in $G$.
\end{corollary}
\begin{proof}
By Theorem \ref{main}, 
$$\mathcal{F}_{\mathcal{A}, \, I_2}(G)=
x_1^n\,\sum_{|\lambda|=n}N_\lambda(G).$$
\end{proof}

The next statement concerns Hamiltonian cycles, or 
connected cycles of length $n=|V(G)|$.

\begin{corollary}
Put $r=n+1$, $x_j=e^{2\pi\sqrt{-1}j/n},\; (j=1, \dots ,n)$ and $t=0$.
Then for any graph $G$ with $n$ vertices
the number of Hamiltonian cycles in $G$ is equal to
$\frac{1}{n}\,\mathcal{F}_{\mathcal{A}, \, I_{n+1}}(G)$.
\label{H}\end{corollary}

\begin{proof}
In this case
$$p_k(x_1,\dots,x_n)=\begin{cases}
0, & k=1, \dots ,n-1,\\ n, & k=n,\end{cases}$$
so that by Theorem \ref{main} 
$$\mathcal{F}_{\mathcal{A}, \, I_{n+1}}(G)=n\,N_{[n]}(G).$$
\end{proof}

\begin{remark}
The value
$\mathcal{F}_{\mathcal{A}, \, I_n+1}(G)$ can be effectively computed for any
graph $G$ as explained in Section \ref{function}. Thus, Corollary \ref{H}
provides a simple algorithm 
for counting the number of Hamiltonian cycles in graphs. Clearly, its 
computational complexity depends on the succession of tensor contractions 
along the edges of $G$. We hope to present a more detailed treatment of this 
problem elsewhere. 
\end{remark}

\end{document}